\documentclass[12pt]{amsart}
\usepackage{epsfig}

\def\N{{\mathbb N}}
\def\Z{{\mathbb Z}}
\def\Q{{\mathbb Q}}
\def\C{\mathbb C}
\def\R{{\mathbb R}}
\def\Hy{{\mathcal H}}
\def\S{{\mathcal S}}

\def\SL{\rm{SL}}
\def\PSL{\rm{PSL}}
\def\tr{\rm{tr}}
\def\Isom{\rm{Isom}}

\title{Singular sets and parameters of generalized triangle orbifolds}
\author{Mikhail Belolipetsky}
\address{Sobolev Institute of Mathematics, Koptyuga 4, 630090 Novosibirsk, Russia \\ and
Max Planck Institute of Mathematics, Vivatsgasse 7, 53111 Bonn, Germany}
\email{mbel@math.nsc.ru}
 
\begin{document}

\begin{abstract}
We study groups generated by three half-turns in the Lobachevsky $3$-space and their
quotient orbifolds. These generalized triangle groups are closely related to 
the arbitrary $2$-generator Kleinian groups. Our main result is a classification
of the singular sets of the generalized triangle orbifolds. We also
present a method to obtain the parameters defining a generalized triangle group
from the structure of the singular set of its quotient orbifold and illustrate it
by examples.
\end{abstract}

\maketitle

\newtheorem{theorem}{Theorem}[section]
\newtheorem{conjecture}[theorem]{Conjecture}
\newtheorem{lemma}[theorem]{Lemma}
\newtheorem{corollary}[theorem]{Corollary}
\theoremstyle{definition}
\newtheorem{definition}[theorem]{Definition}
\newenvironment{remark}{\medskip\noindent\normalsize {\bf Remark.}}{}
\newenvironment{remarkN}[1]{\medskip\noindent\normalsize {\bf Remark #1.}}{}

\section{ Introduction}

After the famous lectures of W.~Thurston \cite{Thur}, the study of
$3$-dimensional manifolds was focused on geometric and, in particular, on hyperbolic
manifolds and orbifolds. An orientable hyperbolic $3$-orbifold can be obtained
as a quotient of the Lobachevsky space $\Hy^3$ by the action of a discrete
subgroup of the group of orientation preserving isometries $\Isom^+\Hy^3= \PSL(2,\C)$. 
This way hyperbolic $3$-orbifolds are related to the discrete subgroups of
$\PSL(2,\C)$ called the {\em Kleinian groups}. Since a subgroup of
$\PSL(2,\C)$ is discrete if and only if each of its $2$-generator subgroups is
discrete (see e.g. \cite{Berdon}), the class of the $2$-generator Kleinian
groups gains a special significance. These groups and, in particular, arithmetic
$2$-generator Kleinian groups were studied by F.~Gehring, C.~Maclachlan,
G.~Martin, J.~Montesinos, A.~Reid and others (see 
\cite{BM, MM3, GM, MM2, ADV, HLM1, J}, \linebreak \cite{MM1}
and the references therein). It appears that the underlying spaces and singular sets of the orbifolds corresponding to the
$2$-generator Kleinian groups may have complicated structure. In a special case of the groups with 
real parameters the orbifolds were studied by J.~Gilman (see \cite{GGM}), 
E.~Klimenko~\cite{K1, K2}; and have been later classified 
by E.~Klimenko and N.~Kopteva (see \cite{KK}).

In this article we consider the groups of hyperbolic isometries generated by
three half-turns in $\Hy^3$. If the axes of the half-turns pairwise intersect
then the group is isomorphic to a Fuchsian triangle group, so the groups
generated by three half-turns can be considered as a generalization of the
Fuchsian triangle groups. These generalized triangle groups are closely related to 
the arbitrary $2$-generator Kleinian groups (see Section~5) but
have much more trackable geometric structure. Thus, since all the three generators
have fixed points (axes) in $\Hy^3$, the fundamental groups of the 
underlying spaces of the corresponding $3$-orbifolds are always trivial, so the
underlying space of such an orbifold is a $3$-sphere $\S^3$ by the Poincar\'e conjecture
(now proved by Perelman). 

Let us briefly describe the contents of the paper.
In Section~2 we give a partially conjectural classification of the singular 
structures of generalized triangle orbifolds. Namely, we show (Theorem~\ref{CT})
that there are eight possible types of singularities and we conjecture that 
under some mild assumptions on the group there are no any other possibilities.
The set of the generalized triangle groups can be parameterized by three
complex numbers that correspond to the complex distances~\cite{Fen} between the
axes of the generators. In Section~3 we use W.~Fenchel's technics
to study this parameterization and its connection with the matrix representation
of the group in $\SL(2,\C)$. Section~4 shows how to find the parameters of the
group for the eight types of singularities from Theorem~\ref{CT}.  This gives an
algorithm which provides a system of three polynomial equations defining
parameters. An explicit connection between generalized triangle groups and 
$2$-generator Kleinian groups with their usual parameterization is given in 
Section~5. In Section~6 we present examples of generalized triangle orbifolds with 
different types of singularities and in Section~7 we consider an exceptional 
example for which Property~RE defined in Section~2 does not hold.
\bigskip

The results of this paper were previously available as a preprint \cite{Bel2},
the current version contains only some changes in the exposition. While working
on the paper I enjoyed the hospitality of MPIM in Bonn. I would like to thank Prof. Alexander 
D.~Mednykh for suggesting this problem for my PhD thesis and for many helpful 
discussions.

The current version of the paper, except for the grammatical corrections and a few small comments,
dates back to February of 2002. It was submitted for a publication in a journal but the 
referee did not find it sufficiently interesting. I also thought that the results were not 
sufficiently strong and decided not to resubmit the paper. The purpose of the present
update is to improve the readability of the text. 
 
\section{ Structure of the singular sets}

The underlying space of a generalized triangle orbifold is the $3$-sphere $\S^3$, so the singular sets 
of the orbifolds are knotted graphs in $\S^3$ with the orbifold 
fundamental groups generated by three involutions. While depicting the singular sets 
we follow a common notation and write the indexes of singularities near the corresponding
components omitting the index $2$. We shall employ the Wirtinger presentation of 
the fundamental group of a knotted graph \cite{Fox}, adapted for the groups of the
orbifolds as in \cite{HQ}, so the elementary arcs of the graph will correspond to the
words in the group generators.

\begin{definition} We call by a {\em generalized triangle group} a discrete subgroup 
$\Gamma$ of $\Isom^+\Hy^3$ generated by three involutions (half-turns). The 
corresponding quotient orbifold ${\mathcal O} = \Isom^+\Hy^3/\Gamma$ is called a 
{\em generalized triangle orbiford.} We say that a generalized triangle orbifold
(and group) has {\em Property~RE} if the rank of the fundamental group of the complement 
of its singular set in $\S^3$ is equal to the rank of~$\Gamma$.
\end{definition}

We shall mainly consider the groups with Property~RE. In some situations this allows us to
use the fundamental group of the complement of the singular set instead of the
orbifold group and saves from certain difficulties which can arise while working with the
groups that have many involutions. An example of an exceptional generalized triangle 
group that does not have Property~RE will be given in the last section. It is
easy to see that the class of groups with Property~RE is closed with respect 
to the Reidemeister moves applied to the singular sets, hence we can safely work 
with the singular graphs remaining inside the class.

\begin{figure}[ht]
\psfig{file=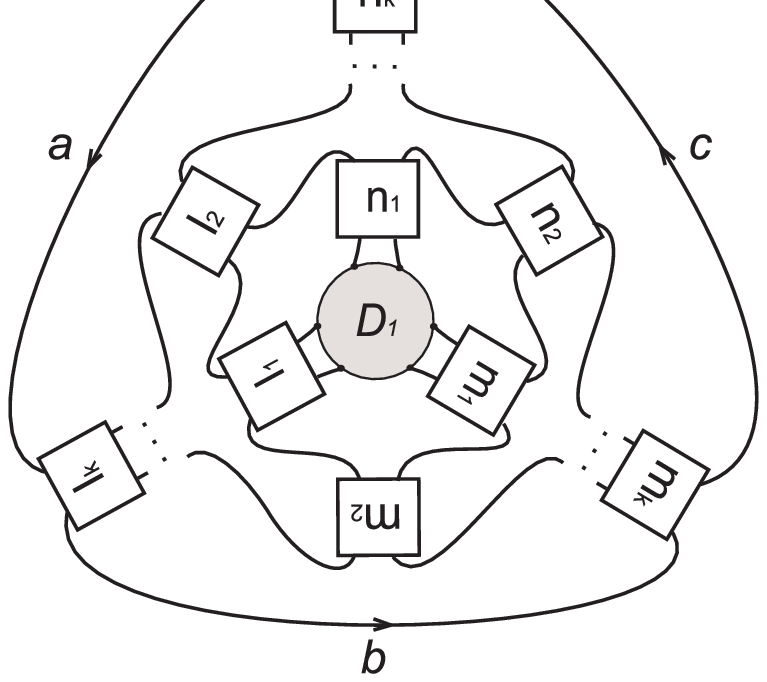, scale=.65} 
\caption{\quad The structure of the singular set of signature
$S(l_1, m_1, n_1, l_2, m_2, n_2,\dots ,l_k, m_k, n_k),$ 
$l_i$, $m_i$, $n_i\in\Z$ for $i=1,\dots,k.$
The boxes denote the integral tangles with the inscribed number of (oriented) crossings.}
\label{SingStr1}
\end{figure}

\begin{figure}[ht]
\psfig{file=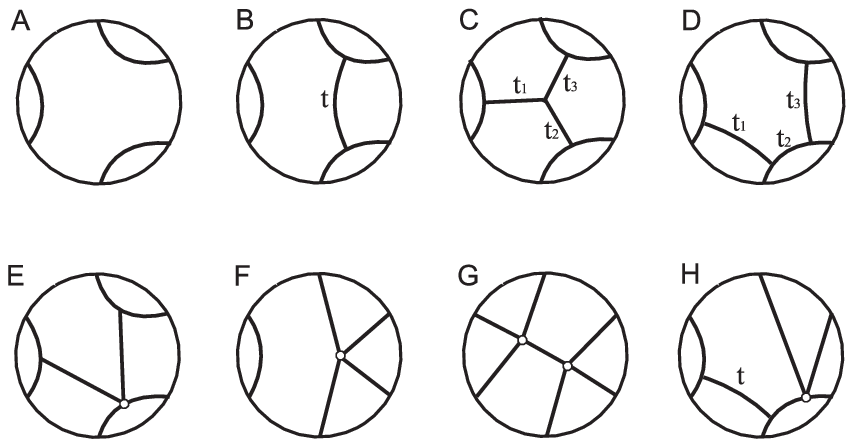, scale=.75} 
\caption{Singular structures inside $D_1$ (here $t_i\in\N\cup\{\infty\}$).}
\label{SingStr2}
\end{figure}

\begin{theorem}\label{CT} (The Classification Theorem)
The singular sets in $\S^3$ that have a stratified structure determined by the signature
$S(l_1, m_1, n_1, l_2,$ $m_2, n_2,\dots ,l_k, m_k, n_k)$ (Figure~\ref{SingStr1}) 
together with the central part $D_1$ of one of the types $A$--$H$ on
Figure~\ref{SingStr2} correspond to generalized triangle orbifolds with
Property~RE.

The orbifolds obtained from the types $A$, $B$, $C\ ($with $1/t_1+1/t_2+1/t_3 > 1)$ and $D\ ($with
$1/t_1+1/t_2 > 1/2$ and $1/t_2+1/t_3 > 1/2)$ are compact; while the orbifolds corresponding to the types
$C\ ($with $1/t_1+1/t_2+1/t_3=1)$, $D\ ($with $1/t_1+1/t_2=1/2$ or
$1/t_2+1/t_3=1/2)$, $E$, $F$, $G$ and $H$ are non-compact.
\end{theorem}

\begin{proof}
The proof is straightforward. We can apply Wirtinger's algorithm to see that in each of the
cases the group has a  presentation with the generators $a,$ $b,$ $c$ corresponding to the
arcs labeled by the same symbols on Figure~\ref{SingStr1}. Applying the same steps
to the complement of the singular set in $\S^3$ we see that the fundamental group of the
complement can be generated by the elements corresponding to the same arcs, so 
Property~RE is satisfied automatically. To distinguish between the compact and non-compact
cases for each of the vertices of the singular graph we can check whether the boundary of the
neighborhood of the vertex is a sphere or a plane, depending on the corresponding subgroup 
of the monodromy group.
\end{proof}

\begin{conjecture} \label{CC}
The classification theorem describes all possible singular sets of the hyperbolic generalized 
triangle orbifolds with Property~RE.
\end{conjecture}

In \cite{Bel2} this conjecture was a part of the classification theorem however
the proof given there is incomplete and contains serious gaps. Let us briefly recall the
argument since it still can be considered as a basis for Conjecture~\ref{CC}.

We begin with a singular graph of a generalized triangle orbifold with Property~RE. Let us
choose three arcs of the graph which correspond to the generators of the group. Using
the Reidemeister moves we can push the arcs outside of the disc which contains the plane projection
of the singular set. This is the first step to obtain the stratified structure. Then we
can show that since we are not allowed to add more generators, the singular arcs can not be
knotted with themselves, and we also can prove that by the Reidemeister moves it is possible
to gather all the vertices of the singular graph in the central part of the diagram. 
This essentially leads to the stratified structure of the singular set depicted on 
Figure~\ref{SingStr1}. It remains to prove that all the tangles on 
Figure~\ref{SingStr1} are integral or, equivalently, that the tangles of the form 
$1/n,$ $n\ge 2$ cannot arise. This may seem obvious from the first view but I do not 
know how to prove this fact. It may require a different approach. 
Finally, it remains to understand the possible singular structures in the center of 
the diagram and here Property~RE is essential, as it is confirmed by the examples from Section~7.

\begin{remarkN}1 Some relatively simple singular structures from
Theorem~\ref{CT} correspond to non-hyperbolic orbifolds (see Examples~6C, 6D). All
such geometric orbifolds can be found in \cite{Dun} and we shall mainly be interested 
in the hyperbolic case.
\end{remarkN}\medskip

\begin{remarkN}2 By the rigidity argument, a generalized triangle orbifold is uniquely 
determined by the signature and the central part from the classification theorem, however, 
different signatures may correspond to the isomorphic orbifolds.
\end{remarkN}

\section{ The matrix representation and the parameters}

In order to obtain a matrix representation of a group generated by
three half-turns one has to fix the axes of the generators, which can be done
by means of the complex distances between the hyperbolic lines
(see \cite[p.~67-70]{Fen}). We define the {\it complex distance} between
two oriented hyperbolic lines as follows: its real part equals the
length of the common perpendicular to the lines and the imaginary part
is given by the angle between the lines, taken from the first to the
second line respecting orientations.

We call three complex distances between the axes of the generators by
the {\it parameters} of the group $\Gamma=\langle a,b,c\rangle$ 
(this denotes a group generated by $a,$ $b,$ $c$ which has relations which are not specified).
Let us note that starting from here by group $\langle a,b,c\rangle$
we mean the group with the given generators taken in a fixed order, so 
we work with the {\em marked groups} generated by $a, b, c$.

We have:
$$
par(\langle a,b,c\rangle)=(\mu (a,b),\ \mu (a,c),\ \mu (c,b))=(\mu_0,\ \mu_1,\ \mu_2),
$$
where $\mu(\cdot,\cdot)=x+iy$, $x\in \R^+$, $y\in [0;2\pi)$ denotes the complex
distance between the oriented axes of the half-turns.

\begin{remark}
Since the directions of the axes of the half-turns are undefined, the parameters are
defined up to a change of the orientation. We can canonically link the
orientation of a hyperbolic line with the matrix of the half-turn in
$\SL(2,\C)$ (see \cite[p.~63]{Fen}). This means that fixing the directions
of the axes of generators is equivalent to choosing the inverse image of $\Gamma < \PSL(2,\C)$ 
under the canonical projection $P:\SL(2,\C)\to\PSL(2,\C)$. Following a usual
agreement we shall call this inverse image by a {\em representation} of our group in $\SL(2,\C)$.
\end{remark}

\begin{lemma}\label{hex}
The complex vector $par(\langle a,b,c\rangle)$ determines the group
$\langle a,b,c\rangle$ uniquely up to a conjugation in $\PSL(2,\C)$, i.e.
if $par(\langle a,b,c\rangle)=par(\langle a',b',c'\rangle)$, then there
exists $h\in\PSL(2,\C)$ such that $a'=hah^{-1}$, $b'=hbh^{-1}$,
$c'=hch^{-1}$.
\end{lemma}

\begin{proof}
Consider the axes of the half-turns $a, b, c$ in $\Hy^3$. We endow the axes by
the orientations and join them in couples by the common perpendiculars. The result
is a hyperbolic right-angled hexagon (see Figure~3).

\begin{figure}[ht]
\psfig{file=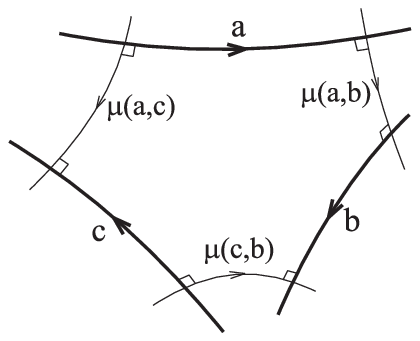, scale=.75}
\caption{The axes of the half-turns $a, b, c$ together with their common
perpendiculars define a right-angled hyperbolic hexagon.}
\end{figure}

The vector $par(\langle a,b,c\rangle)$ represents the three complex lengths of the
sides of our hexagon which are common perpendiculars of the given axes. By
a theorem from \cite[p.~94]{Fen} three pairwise non-adjacent sides define
right-angled hexagon uniquely up to a simultaneous change of orientations
of the other three sides. It means that any other hexagon defined by
another group with the same parameters can be translated to the given
one by a hyperbolic isometry. This isometry obviously defines the required
conjugation $h\in\PSL(2,\C)$.

It may seem that we have missed the possibility of changing the
orientations of the axes but we have already noted that the orientations
affect only representation of the group in $\SL(2,\C)$ but not the isometries
themselves. So in this situation we can change the orientations as it is
needed.
\end{proof}

The next step is to find some convenient representation of the group 
$\langle a,b,c\rangle$ in $\SL(2,\C)$ and describe it in terms of
$par(\langle a,b,c\rangle)$. We denote the matrices corresponding to
the half-turns about $a$, $b$, $c$ by the capitals $A$, $B$, $C$.

A half-turn about an oriented hyperbolic line $m$ is represented by the matrix
$M\in\SL(2,\C)$ of the form $(\begin{array}{cc} m_{11}&m_{12}\\m_{21}&-m_{11}\end{array})$,
which is characterized by $M^2=-I$ or, equivalently, $\tr(M)=0$.
This matrix $M$ is called {\it the normalized line-matrix} of the line $m$.
By extending the transformation $M$ to the boundary $\overline{\C}=\partial\Hy^3$,
it is easy to find the set of its fixed points $\it{fix}(M) \subset \overline{\C}$. 
Geometrically, the fixed points $\{z_1,\ z_2\}=\it{fix}(M)$ are the
ends of the arc orthogonal to $\C$ which represents the line $m$ in the
upper half-space model of $\Hy^3$.

After a suitable conjugation of the group $\langle A,B,C\rangle$ in $\SL(2,\C)$
one can suppose that $\it{fix}(A)=\{ 1/\beta ;-1/\beta\}$ and 
$\it{fix}(B)=\{ \beta ;-\beta\}$ with $\beta\in\C$, $|\beta|\geq 1$.

We have:
$$
A=\left(\begin{array}{cc} 0&i/\beta\\i\beta&0\end{array}\right),\ \
B=\left(\begin{array}{cc} 0&i\beta\\i/\beta&0 \end{array}\right),\ \
C=\left(\begin{array}{cc} c_{11}&c_{12}\\c_{21}&-c_{11}\end{array}\right),
$$
where $-c_{11}^2+c_{12}c_{21}=1$.

To find $\beta$ and $c_{ij}$ in terms of $par(\langle a,b,c\rangle)$ we
use the following basic formula allowing us to express the complex
distance between the hyperbolic lines $m_1$, $m_2$ in terms of the matrices
$M_1$, $M_2$ of the lines (see \cite[p.~68]{Fen}):
$$\cosh(\mu(m_1,m_2))=-\frac{1}{2} \tr\,(M_1M_2).$$

Denote $\rho_k=-2\cosh(\mu_k)$, $k=0$, $1$, $2$. We have
$$\rho_0=\tr(AB)=\tr\left(\begin{array}{cc} -1/\beta^2&0\\0&-\beta^2\end{array}\right);$$
$$\rho_1=\tr(AC)=\tr\left(\begin{array}{cc} ic_{21}/\beta&-ic_{11}/\beta\\ic_{11}\beta&ic_{12}\beta\end{array}\right);$$
$$\rho_2=\tr(BC)=\tr\left(\begin{array}{cc} ic_{21}\beta&-ic_{11}\beta\\ic_{11}/\beta&ic_{12}/\beta\end{array}\right).$$

Hence we obtain a system of equations on $\beta$ and $c_{ij}$:
$$
\left\{
\begin{array}{l}
{\displaystyle -1/\beta^2-\beta^2=\rho_0};\\
{\displaystyle ic_{21}/\beta+ic_{12}\beta=\rho_1};\\
{\displaystyle ic_{21}\beta+ic_{12}/\beta=\rho_2};\\
{\displaystyle c_{11}=i\sqrt{c_{12}c_{21}+1}}.
\end{array}
\right.
$$
Solving these equations we can find all the parameters in our representation.
Let us note that non-uniqueness of the solution of the system in general does not
preoccupy us because different solutions will give different representations of the group
$\Gamma=\langle a,b,c\rangle$ which are conjugate in $\SL(2,\C)$.
The latter statement is a direct consequence of Lemma~\ref{hex}, since by the
construction different solutions correspond to the same $par(\Gamma)$.
For the sake of concreteness we can always fix the analytic branches of the square
roots in the following formulas.

We now can write down the representation of $\langle a,b,c\rangle$
in $\SL(2,\C)$:

$$
A=\left(\begin{array}{cc} 0&i/\beta\\i\beta&0\end{array}\right),\ \
B=\left(\begin{array}{cc} 0&i\beta\\i/\beta&0 \end{array}\right),\ \
C=\left(\begin{array}{cc} c_{11}&c_{12}\\c_{21}&-c_{11}\end{array}\right);
$$
\begin{equation}
\begin{array}{ll}
\phantom{wwww}&{\displaystyle
  \beta=\sqrt{\frac{-\rho_0+\sqrt{\rho_0^2-4}}{2}}},\\
\phantom{wwww}&{\displaystyle \vphantom{\sqrt{\frac{-\rho_0+\sqrt{\rho_0^2-4}}{2}}}
  c_{21}=\frac{\rho_1/\beta-\rho_2\beta}{i/\beta^2-i\beta^2},\ \
  c_{12}=\frac{-\rho_1\beta+\rho_2/\beta}{i/\beta^2-i\beta^2},\ \
  c_{11}=i\sqrt{c_{12}c_{21}+1};}\\
\phantom{wwww}&{\displaystyle \vphantom{\sqrt{\frac{-\rho_0+\sqrt{\rho_0^2-4}}{2}}}
  \rho_k=-2\cosh(\mu_k)\ {\rm for}\ k=0,1,2}.
\end{array}
\end{equation}

To obtain the representation of $\langle a,b,c\rangle$ in $\PSL(2,\C)$ one has
to consider the image of the $\SL(2,\C)$-representation under the
canonical projection $P:\SL(2,\C)\to\PSL(2,\C)$.

\section{ Deducing parameters from the singular structure}

Knowing the structure of the singular set of a generalized triangle orbifold
we can obtain a presentation of its fundamental group and a representation of
the group in $\SL(2,\C)$. This representation is defined by three complex parameters
$\rho_0$, $\rho_1$ and $\rho_2$, so it is possible to find the parameters
corresponding to a given singular structure. In this section and in
Section~6 we shall show how this can be carried out for the generalized triangle groups 
with Property~RE.

Suppose we are given a singular set of a generalized triangle orbifold corresponding
to one of the cases of Theorem~\ref{CT}.
Using Wirtinger's algorithm we can start from the given three arcs
$a$, $b$, $c$ and passing through the tangles obtain words in $a$, $b$, $c$
corresponding to all the other arcs. It is clear that until we reach the
central part of the diagram we do not obtain any equations on the
generators, so all the words in the presentation of the group come from
the central part of the diagram. Equations defining parameters are 
then obtained by taking traces of corresponding elements in $\SL(2,\C)$. With the help 
of the well known formulas
\begin{equation}
\begin{array}{l}
\tr(W_1 W_2) = \tr(W_2 W_1),\\
\tr(W_1 W_2 W_1^{-1}) = \tr(W_2),\\
\tr(W_1 W_2)=\tr(W_1)\tr(W_2)-tr(W_1^{-1} W_2),\\
\tr(H^{-1})=-\tr(H) \ \ \ {\rm if\ }H{\rm\ is\ a\ matrix\ of\ a\ half-turn,}\\
\end{array}
\end{equation}
the traces of the expressions in $A$, $B$, $C$ and their inverses can be written in terms of
$\tr(AB)=\rho_0$, $\tr(AC)=\rho_1$, $\tr(BC)=\rho_2$ and $\tr(ABC)$ (an
explicit expression for $\tr(ABC)$ in terms of $\rho_i$ can be obtained using
(1)). Since we have three complex parameters, we require three independent
equations to define them. In the remaining part of this section we shall give
the presentations of the groups and corresponding equations for each of the
types $A$--$H$\ from Theorem~\ref{CT}.

\bigskip
\psfig{file=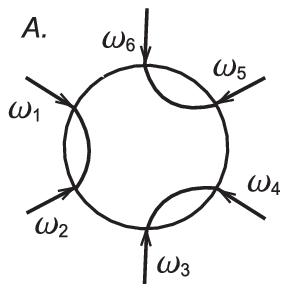, scale=.75}
\bigskip
\hfill
\parbox[b]{0.7\textwidth}{
Group:\\
$\langle a,\ b,\ c \mid a^2,\ b^2,\ c^2,\ w_1 w_2,\ w_3 w_4 \rangle,$

\medskip
equations:

\smallskip
$\left\{
\begin{array}{l}
\tr(W_1 W_2)=2,\\ \tr(W_3 W_4)=2,\\ \tr(W_5 W_6)=2
\end{array}
\right.$

\smallskip
(We denote by $w_i$ the words in $a$, $b$, $c$\ and by $W_i$~--- corresponding
products in $\SL(2,\C)$, $i=1,\dots,6.$ )
}

It follows from \cite{Fox} that one of the equations $w_1 w_2$, $w_3 w_4$,
$w_5 w_6$ in the fundamental group is always a corollary of the two others.
However, by the rigidity argument in the matrix group all the three words
are independent and give three independent equations on the parameters
(see Example~6A for an illustration).

\bigskip
\psfig{file=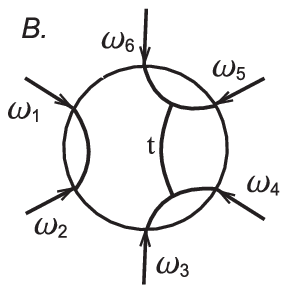, scale=.75}
\hfill
\parbox[b]{0.7\textwidth}{
Group: \\
$\langle a,\ b,\ c \mid a^2,\ b^2,\ c^2,\ w_1 w_2,\ (w_3 w_4)^t,\ (w_5 w_6)^t\rangle,$

\medskip
equations:

\smallskip
$\left\{
\begin{array}{l}
\tr(W_1 W_2)=2,\\ \tr(W_3 W_4)=-2cos(\pi/t),\\ \tr(W_5 W_6)=-2cos(\pi/t).
\end{array}
\right.$
\bigskip
}

\bigskip
\psfig{file=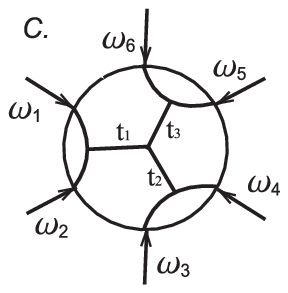, scale=.75}
\hfill
\parbox[b]{0.7\textwidth}{
Group: \\
$\langle a,\ b,\ c \mid a^2,\ b^2,\ c^2,\ (w_1 w_2)^{t_1},\ (w_3 w_4)^{t_2},\ (w_5 w_6)^{t_3}\rangle,$

\medskip
equations:

\smallskip
$\left\{
\begin{array}{l}
\tr(W_1 W_2)=-2cos(\pi/t_1),\\ \tr(W_3 W_4)=-2cos(\pi/t_2),\\ \tr(W_5 W_6)=-2cos(\pi/t_3).
\end{array}
\right.$
\bigskip
}

\bigskip
\psfig{file=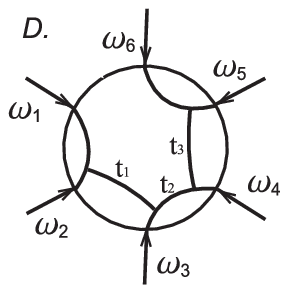, scale=.75}
\hfill
\parbox[b]{0.7\textwidth}{
Group: \\
$\langle a,\ b,\ c \mid a^2,\ b^2,\ c^2,\ (w_1 w_2)^{t_1},\ (w_5 w_6)^{t_2},\ (w_1 w_2 w_3)^{t_3}\rangle,$

\medskip
equations:

\smallskip
$\left\{
\begin{array}{l}
\tr(W_1 W_2)=-2cos(\pi/t_1),\\ \tr(W_5 W_6)=-2cos(\pi/t_3),\\
\tr(W_1 W_2 W_3)=-2cos(\pi/t_2).
\end{array}
\right.$
\bigskip
}

It can be seen that the word $W_1 W_2 W_3$ has an odd length, which implies that in
this case we require the expression for $\tr(ABC)$. In fact, this is the only type for which
we may need this expression.

\bigskip
\psfig{file=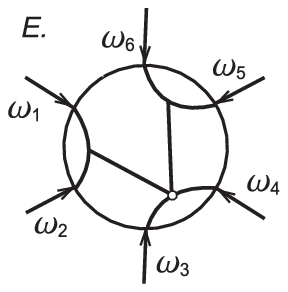, scale=.75}
\hfill
\parbox[b]{0.7\textwidth}{
Group: \\
$\langle a,\ b,\ c \mid a^2,\ b^2,\ c^2,\ (w_1 w_2)^2,\ (w_5 w_6)^2\rangle,$

\medskip
equations:

\smallskip
$\left\{
\begin{array}{l}
\tr(W_1 W_2)=0,\\ \tr(W_5 W_6)=0,\\ \tr(W_3 W_4)=-2.
\end{array}
\right.$
\bigskip
}

The third equation follows from the fact that the multiple of two
half-turns with parallel axes is a parabolic isometry, so $W_3 W_4$ is
a parabolic corresponding to the cusp. We choose the sign $'-'$ for the trace
of the parabolic element because the axes of the half-turns have the same
directions.

\bigskip
\psfig{file=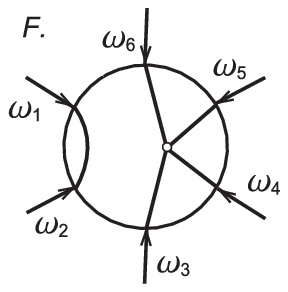, scale=.75}
\hfill
\parbox[b]{0.7\textwidth}{
Group: \\
$\langle a,\ b,\ c \mid a^2,\ b^2,\ c^2,\ w_1 w_2\rangle,$

\medskip
equations:

\smallskip
$\left\{
\begin{array}{l}
\tr(W_1 W_2)=2,\\ \tr(W_3 W_4)=-2,\\ \tr(W_5 W_6)=-2.
\end{array}
\right.$
\bigskip
}

\bigskip
\psfig{file=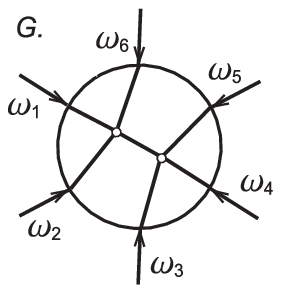, scale=.75}
\hfill
\parbox[b]{0.7\textwidth}{
Group: \\
$\langle a,\ b,\ c \mid a^2,\ b^2,\ c^2,\ (w_3 w_4 w_5)^2 \rangle,$

\medskip
equations:

\smallskip
$\left\{
\begin{array}{l}
\tr(W_3 W_4 W_5)=0,\\ \tr(W_1 W_2)=-2,\\ \tr(W_3 W_4)=-2.
\end{array}
\right.$
\bigskip
}

Here, as in the case $D$, we again have a word $W_3 W_4 W_5$ of an odd length, but
since $\tr(W_3 W_4 W_5)=0$, in this case we can easily avoid using the explicit expression
for $\tr(ABC)$. Indeed, using formulas $(2)$ we can present
$\tr(W_3 W_4 W_5)$ as $\tr(ABC){\rm P}[\tr(AB), \tr(AC), \tr(BC)]$, where
${\rm P[...]}$ denotes a polynomial, and since $\tr(ABC)\neq 0$ we obtain the equation
$\rm P[\tr(AB), \tr(AC), \tr(BC)]=0$. The same procedure can be applied for the type $H$
as well.

\bigskip
\psfig{file=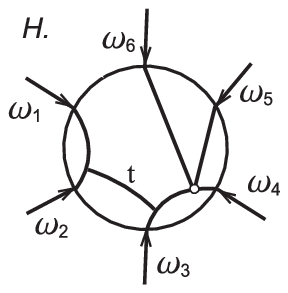, scale=.75}
\hfill
\parbox[b]{0.7\textwidth}{
Group: \\
$\langle a,\ b,\ c \mid a^2,\ b^2,\ c^2,\ (w_1 w_2)^t,\ (w_1 w_2 w_3)^2 \rangle,$

\medskip
equations:

\smallskip
$\left\{
\begin{array}{l}
\tr(W_1 W_2)=-2cos(\pi/t),\\ \tr(W_5 W_6)=-2,\\
\tr(W_1 W_2 W_3)=0.
\end{array}
\right.$
\bigskip
}

\section{ Relation to the two-generator Kleinian groups}

Let $\Gamma_0=\langle f,g\rangle$ be a discrete subgroup of $\PSL(2,\C)$.
There is a natural way to associate to $\Gamma_0$ a generalized triangle
group $\Gamma$ which is equal to $\Gamma_0$ or contains it as a subgroup
of index $2$. 
The latter is controlled by the first homology group of $\Gamma$ having rank
$2$ or $3$, and it can be seen that both cases are possible.
It follows that such properties as discreteness and arithmeticity
of two-generator subgroups of $\PSL(2,\C)$ can be studied
using the generalized triangle groups.

The group $\Gamma$ is constructed as follows (cf. \cite{Gilm}):
If neither $f$ nor $g$ are parabolic then let $N$ be a common perpendicular
to the axes of $f$ and $g$. In the case of a parabolic generator, the corresponding
end of the line $N$ is the fixed point of the parabolic. Clearly,
this construction uniquely determines the line $N$ for any two-generator Kleinian group
$\Gamma_0=\langle f,g\rangle$. Let $c$ be the half-turn about $N$. Then
there exist such half-turns $a$ and $b$ that $f=ac$ and $g=cb$
(by \cite[p.~47]{Fen}). Therefore, we obtain a group $\Gamma=\langle a,b,c\rangle$
which contains $\Gamma_0$ as a subgroup of index at most $2$.

Reciprocally, starting with the group $\Gamma$ generated by three half-turns
$a$, $b$ and $c$ one can easily obtain its two generator subgroup
$\Gamma_0=\langle ab,ac\rangle$. Since $\Gamma=\Gamma_0\cup a\Gamma_0\
(a^2=id)$, the index of $\Gamma_0$ in $\Gamma$ is equal to $2$ or $1$
depending on whether or not the element $a$ is contained in $\Gamma_0$.

The set of two-generator Kleinian groups is usually parameterized by three
complex numbers \cite{GM}: 
$$ \beta(f)=\tr^2(f)-4,\ \beta(g)=\tr^2(g)-4,\ \gamma(f,g)=\tr([f,g])-2. $$ 
A standard computation with traces shows that this
parameters are related to our parameterization for the corresponding
generalized triangle group by the following formulas: 
$$ \beta(f)=\rho_1^2-4,\
\beta(g)=\rho_2^2-4,\
\gamma(f,g)=\rho_0^2+\rho_1^2+\rho_2^2+\rho_0\rho_1\rho_2-4. $$

Since we have described the singular sets of the generalized triangle
orbifolds, we can at least in theory present an algorithm which using
equations from the previous section enumerates the parameters
corresponding to the generalized triangle groups and discrete two-generator
subgroups of $\PSL(2,\C)$. The main practical difficulty is that the
complexity of the equations grows as the singular structures
become more complicated. The other difficulty is that we have to cast away
the redundant solutions of the equations which correspond to non-discrete
subgroups of $\PSL(2,\C)$. In the next section we shall see how this can be 
done for the examples considered there. In general, if an orbifold is arithmetic 
the equations always have only one complex solution (up to a conjugation)
and it defines the parameters, but in the non-arithmetic case there is no any 
general method to choose the right complex place.

Let us also note that the suggested algorithm shows only that our set of
parameters is \emph{semi-recursive}: we can enumerate the parameters of
the discrete generalized triangle groups with Property~RE, but there is no known procedure to decide
in a finite number of steps whether or not a given triple of complex numbers defines
such a group.

\section{ Examples of the generalized triangle orbifolds. }

We now consider examples for different types of singularities from Theorem~\ref{CT}
and find the parameters defining corresponding generalized triangle groups. We
shall try to choose the examples which may deserve a particular interest or previously 
appeared in a different context. Some othere examples of generalized 
triangle orbifolds can be found in \cite{Bel_m} and \cite{Bel1}.

\bigskip

\psfig{file=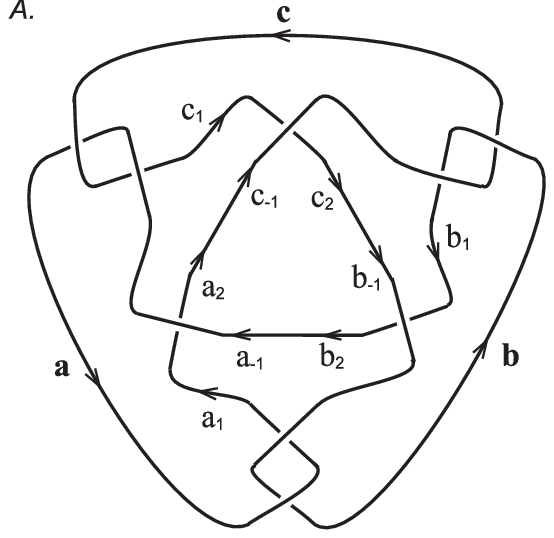, scale=.75}
\hfill
\parbox[b]{0.5\textwidth}{
$\begin{array}{l}
a_{-1} = cac^{-1}\\
b_{-1} = aba^{-1}\\
c_{-1} = bcb^{-1}\\
a_1 = b_{-1}^{-1}ab_{-1}\\
b_1 = c_{-1}^{-1}bc_{-1}\\
b_2 = b_{-1}^{-1}b_1 b_{-1}\\
a_{-1}^{-1}b_2 = ca^{-1}c^{-1}ab^{-1}a^{-1}bc^{-1}bcb^{-1}aba^{-1}\\
\end{array}
$
\bigskip

Matrix group:
\medskip

$\langle A,\ B,\ C \mid -A^2, -B^2, -C^2,\\
CA^{-1}C^{-1}AB^{-1}A^{-1}BC^{-1}BCB^{-1}ABA^{-1},\\
AB^{-1}A^{-1}BC^{-1}B^{-1}CA^{-1}CAC^{-1}BCA^{-1},\\
BC^{-1}B^{-1}CA^{-1}C^{-1}AB^{-1}ABA^{-1}CAB^{-1}
\rangle$
}
\bigskip

Let $\tr(AB)=t_0$, $\tr(AC)=t_1$, $\tr(BC)=t_2$.

By taking traces of the relations we obtain a system of equations on $t_0$, $t_1$
and $t_2$\footnote{We used the computer program $<$tracer$>$(C++, 242
lines) to calculate the trace polynomials for the generalized triangle groups.}:

$$\left\{
\begin{array}{l}
3t_0-t_0^3-2t_0t_1^2+t_0^3t_1^2+2t_0^2t_1t_2-t_0^2t_1^3t_2-2t_0t_2^2+t_0^3t_2^2-t_0^3t_1^2t_2^2-t_0^2t_1t_2^3
= 2,\\
3t_2-t_2^3-2t_2t_0^2+t_2^3t_0^2+2t_2^2t_0t_1-t_2^2t_0^3t_1-2t_2t_1^2+t_2^3t_1^2-t_2^3t_0^2t_1^2-t_2^2t_0t_1^3
= 2,\\
3t_1-t_1^3-2t_1t_2^2+t_1^3t_2^2+2t_1^2t_2t_0-t_1^2t_2^3t_0-2t_1t_0^2+t_1^3t_0^2-t_1^3t_2^2t_0^2-t_1^2t_2t_0^3
= 2.
\end{array}
\right.$$

We are interested only in complex solutions of this system and by the
symmetry of the singular graph all the three roots should be equal.
The only such solutions are (approximately):
$$ (t_0, t_1, t_2) = (0.662359 \pm 0.56228i, 0.662359 \pm 0.56228i,
0.662359 \pm 0.56228i).$$
So we have $\rho_0=\rho_1=\rho_2 = 0.662359... \pm 0.56228... i$.
The signs of the imaginary parts correspond to the orientations of
the generating arcs and do not affect the presentation of the group
in $\PSL(2,\C)$.

The primitive polynomial of the root is $t^3-t+1$. Using the
arithmeticity test for generalized triangle groups from \cite{Bel1}
it can be easily verified that the orbifold is arithmetic. It is also
interesting to note (see \cite{MV1}) that the two-fold covering of this
orbifold is the Weeks-Mateveev-Fomenko manifold~--- the hyperbolic 
$3$-manifold of the smallest volume.

(Here we could also have used the symmetry of the singular graph
from the beginning by initially assuming that $\tr(AB)=\tr(AC)=\tr(BC)$. This
would considerably simplify the calculations but we are mainly interested
in demonstrating the general method.)

\bigskip

\psfig{file=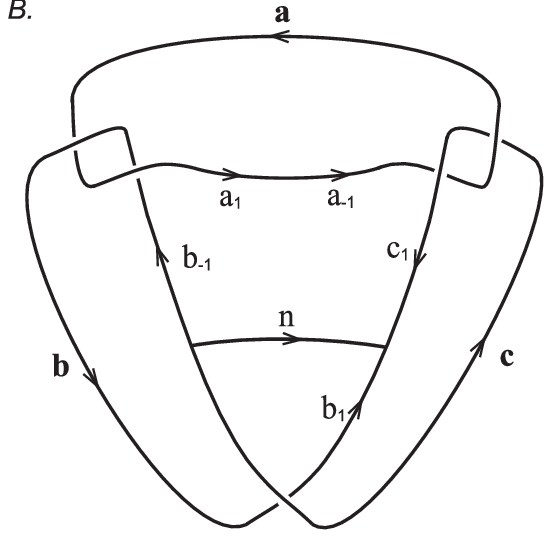, scale=.75}
\hfill
\parbox[b]{0.4\textwidth}{
$\begin{array}{l}
a_1 = bab^{-1}\\
b_1 = cbc^{-1}\\
c_1 = aca^{-1}\\
a_{-1} = c_1^{-1}ac_1\\
b_{-1} = a_1^{-1}ba_1\\
\\
a_{-1}^{-1}a_1 = ac^{-1}a^{-1}ca^{-1}bab^{-1}\\
b_{-1}^{-1}c^{-1} = ba^{-1}b^{-1}ab^{-1}c^{-1}\\
b_1c_1 = cbc^{-1}aca^{-1}\\
\end{array}
$
\vspace{50pt}
}
\bigskip

\noindent Equations on $t_0=\tr(AB)$, $t_1=\tr(AC)$, $t_2=\tr(BC)$:
$$\left\{
\begin{array}{l}
t_0^2t_1^2+t_0t_1t_2-t_0^2-t_1^2+2 = 2,\\
t_0^2t_2+t_0t_1-t_2 = -2\cos(\pi/n),\\
t_1^2t_2+t_0t_1-t_2 = -2\cos(\pi/n).\\
\end{array}
\right.$$

We shall consider only $n=3$. Up to the symmetry and change of orientations, the
system has only one complex solution which defines the complex place of the
orbifold group:
$$\rho_0 = \frac{i+\sqrt{3}}{2},\ \rho_1 = \frac{i-\sqrt{3}}{2},\ \rho_2=0.$$
The squares of $\rho_0$, $\rho_1$ are algebraic integers and by the arithmeticity
test \cite{Bel1} the group is arithmetic and defined over the field $\Q[\sqrt{-3}]$.

\noindent (For $n=2$ and $n=\infty$ the group has no any complex places and the
orbifold is not hyperbolic. For large values of $n$ the group has more than one
complex place induced by the Galois automorphisms $\cos(\pi/n)\to \cos(k\pi/n)$, 
$(k,n)=1$, and so it is not arithmetic.)

\bigskip

\psfig{file=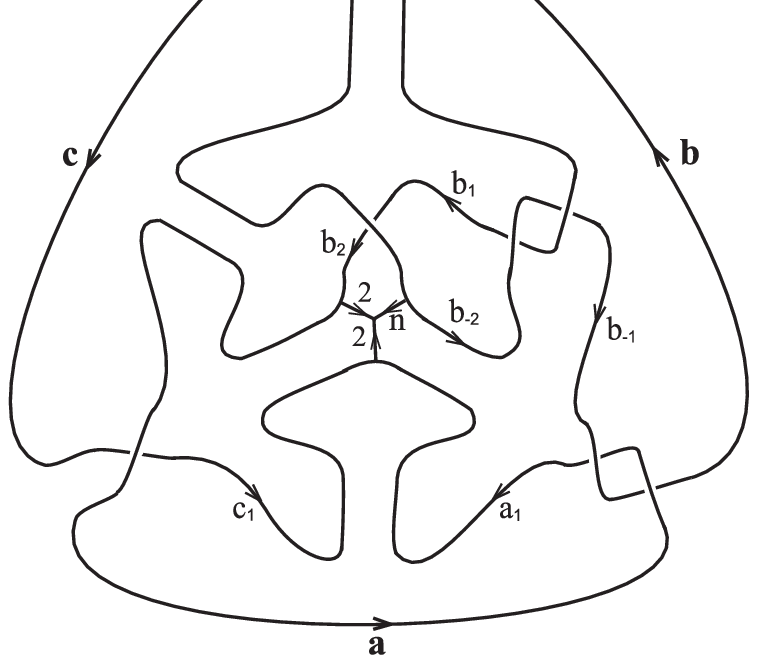, scale=.70}
\hfill
\parbox[b]{0.5\textwidth}{
$\begin{array}{l}
b_{-2} = b^{-1}aba^{-1}b\\
b_{-1} = aba^{-1}\\
b_1 = b_{-2}bb_{-2}^{-1}\\
b_2 = cb_1c^{-1}\\
a_1 = b_{-1}^{-1}ab_{-1}\\
c_1 = aca^{-1}\\
\\
c_1a_1 = acb^{-1}aba^{-1}\\
b_{-2}^{-1}c^{-1} = b^{-1}ab^{-1}a^{-1}bc^{-1}\\
b_2a^{-1} = \\
\ \ cb^{-1}aba^{-1}bab^{-1}a^{-1}bc^{-1}a^{-1}
\end{array}
$
\vspace{90pt}
}
\bigskip

\noindent Equations on $t_0=\tr(AB)$, $t_1=\tr(AC)$, $t_2=\tr(BC)$:
$$\left\{
\begin{array}{l}
t_0t_2+t_1 = 0,\\
t_0^2t_2+t_0t_1-t_2 = -2\cos(\pi/n),\\
(t_0^3-3t_0)(t_0t_1t_2+t_1^2+t_0^2-2)+t_0t_1^2+t_1t_2-t_0 = 0.
\end{array}
\right.$$
The complex solution ($n>3$):
$$\begin{array}{l}
\rho_0^2 = \frac{1}{2}(1+4\sin^2(\pi/n)) \left(1\pm\sqrt{1-\frac{4}{1+4\sin^2(\pi/n)}}\right),\\
\rho_1 = -2\cos(\pi/n)\rho_0,\\
\rho_2 = 2\cos(\pi/n).
\end{array}
$$
(Here $\rho_2$ is always a real number; for $n=2,3$, $\rho_0$ and $\rho_1$ are
also real, so the corresponding orbifolds are not hyperbolic.
According to \cite{Dun}, for $n=2$ the orbifold admits the spherical structure 
and for $n=3$ it is Euclidean.)

As it was shown in \cite{RV}, the orbifolds considered here can be obtained
as the quotient orbifolds of the Fibonacci manifolds $F(2n,n)$ by their full
groups of isometries. Hence using the arithmeticity test from \cite{Bel1}
we can now easily find all arithmetic Fibonacci manifolds (see also
\cite{HLM}). The field of definition of the orbifold group is $\Q[\rho_0^2,
\cos^2(\pi/n)]$. It has exactly one complex place only if $n=4$, $5$, $6$, $8$
$12$, $\infty $; and it is a matter of a direct verification that these cases
satisfy the remaining conditions of the arithmeticity test. Thus,
the orbifolds and corresponding Fibonacci manifolds are arithmetic
for $n=4$, $5$, $6$, $8$, $12$, $\infty $.

\bigskip

\psfig{file=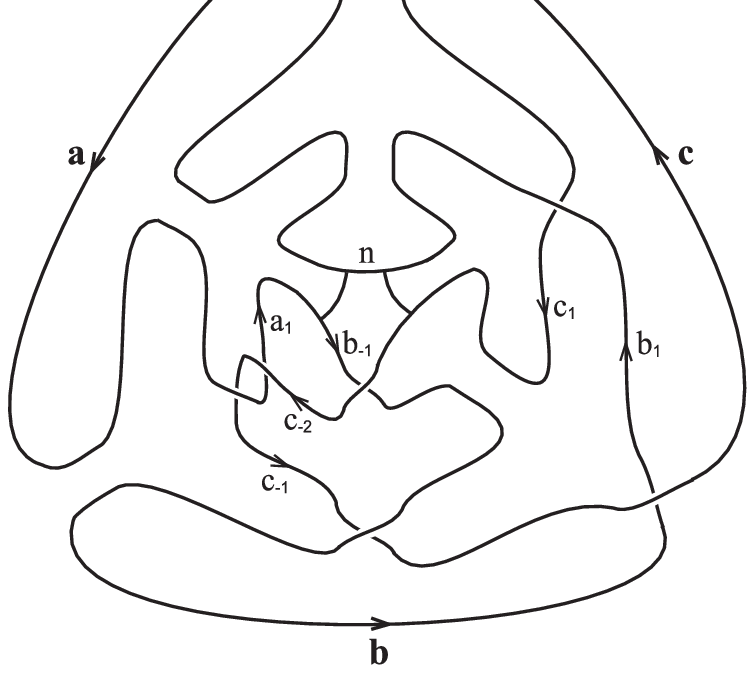, scale=.70}
\hfill
\parbox[b]{0.5\textwidth}{
$\begin{array}{l}
c_{-2} = ac_{-1}a^{-1}\\
c_{-1} = b^{-1}cb\\
b_{-1} = c_{-2}^{-1}bc_{-2}\\
a_1 = c_{-2}^{-1}ac_{-2}\\
b_1 = cbc^{-1}\\
c_1 = b_{1}cb_{1}^{-1}\\
\\
a_1b_{-1}^{-1} = ab^{-1}c^{-1}ab^{-1}cba^{-1}\\
c_{-2}^{-1}c_1 = ab^{-1}c^{-1}ba^{-1}cbcb^{-1}c^{-1}\\
a^{-1}a_1b_{-1}^{-1} = b^{-1}c^{-1}ab^{-1}cba^{-1}
\end{array}
$
\vspace{110pt}
}
\bigskip

\noindent Equations on $t_0=\tr(AB)$, $t_1=\tr(AC)$, $t_2=\tr(BC)$:
$$\left\{
\begin{array}{l}
t_0 = 0,\\
4t_2^2-t_1^2t_2^2-t_2^4+t_1^2-2 = 0,\\
4t_1^2-t_1^2t_2^2-t_1^4 = 4cos^2(\pi/n).\\
\end{array}
\right.$$
The equation $\tr(AB) = 0$ geometrically means that the axes of the half-turns
$a$ and $b$ intersect orthogonaly in $\Hy^3$.

The system has exactly one pair of complex-conjugate roots for any
$n>2$ (for $n=2$ the orbifold is Euclidean), these roots give the
values of the parameters $\rho_0 (=0),$ $\rho_1,$ $\rho_2$ defining the
orbifold group. Since the analytic formulas for the solutions are
rather complicated we do not present them here. One can see that the
field of definition $\Q[\rho_1^2, \rho_2^2]$ of the group  has
exactly one complex place iff $n=3$, $4$, $6$, and it is easy to
check the other conditions of the arithmeticity test \cite{Bel1}
for these values of $n$. So we obtain that the orbifolds are arithmetic
for $n=3$, $4$, $6$.

These groups and orbifolds first appeared in \cite{HKM2} as an infinite
one-parameter family extending one of the three regular tessellations of
$\Hy^3$. It is interesting to note that the corresponding family for the one of
the two remaining regular tessellations is also well known \cite{HKM1} and consists
of the Fibonacci manifolds considered in the previous example. For
the third tessellation the question of constructing such a family is still
open.

\bigskip

\psfig{file=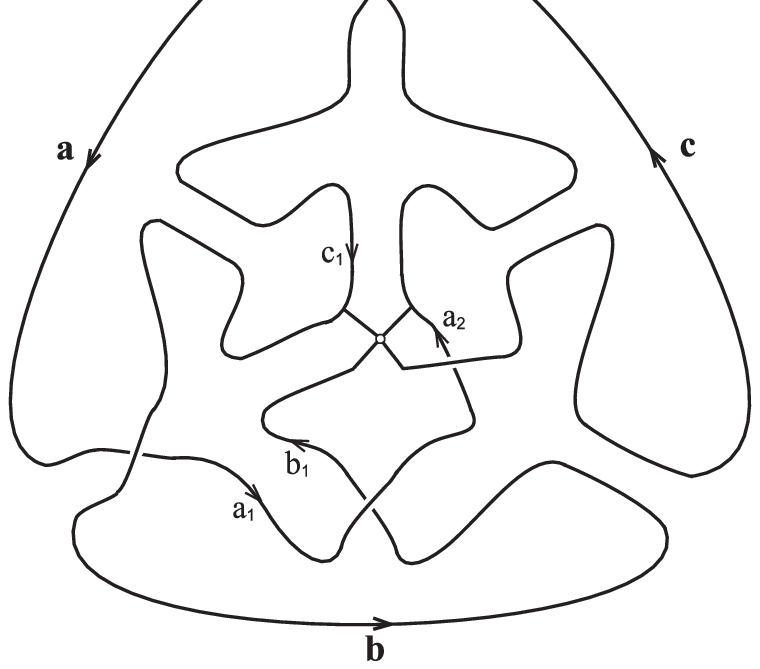, scale=.70}
\hfill
\parbox[b]{0.5\textwidth}{
$\begin{array}{l}
a_1 = bab^{-1} \\
b_1 = a_1ba_1^{-1} \\
c_1 = aca^{-1} \\
a_2 = ca_1c^{-1} \\
\\
c_1b^{-1} = aca^{-1}b^{-1} \\
a_2a^{-1} = cbab^{-1}c^{-1}a^{-1} \\
b_1c^{-1} = baba^{-1}b^{-1}c^{-1} \\
\end{array}
$
\vspace{136pt}
}
\bigskip

\noindent Equations on $t_0=\tr(AB)$, $t_1=\tr(AC)$, $t_2=\tr(BC)$:
$$\left\{
\begin{array}{l}
t_0t_1+t_2 = 0,\\
t_0t_1t_2+t_0^2+t_1^2-2 = 0,\\
-t_0^2t_2-t_0t_1+t_2 = -2.
\end{array}
\right.$$
The solution:
$$\rho_0^2 = \rho_1^2 = 1+i,\ \rho_2^2=2i.$$

This orbifold first appeared in \cite{Adams} as the third hyperbolic
$3$-orbifold with a non-rigid cusp. We see that the orbifold is
arithmetic with the field of definition $\Q[i]$ and so its group is
commensurable in $\PSL(2,\C)$ with the Picard group $\PSL(2,\Z[i])$.

\bigskip

\psfig{file=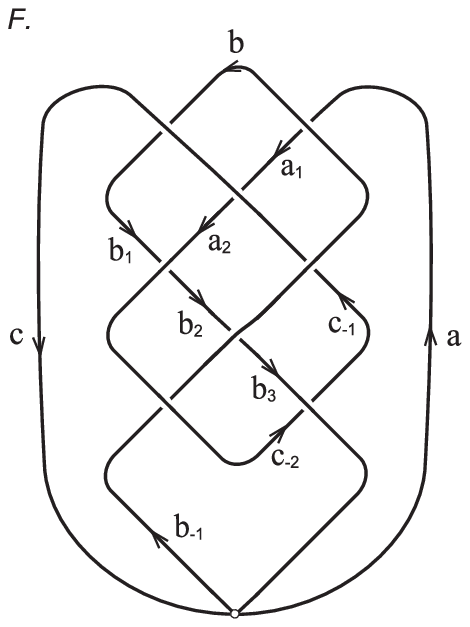, scale=.75}
\hfill
\parbox[b]{0.5\textwidth}{
$\begin{array}{l}
a_1 = bab^{-1}\\
a_2 = ca_1c^{-1}\\
b_1 = cbc^{-1}\\
b_2 = a_2b_1a_2^{-1}\\
b_3 = b^{-1}b_2b\\
c_{-1} = b^{-1}cb\\
c_{-2} = b_3^{-1}c_{-1}b_3\\
b_{-1} = a_2^{-1}ba_2\\
\\
c_{-2}^{-1}a_2 = \\
\  b^{-1}cbab^{-1}a^{-1}b^{-1}c^{-1}baba^{-1}b^{-1}c^{-1}bcbab^{-1}c^{-1}\\
b_{-1}^{-1}c = cba^{-1}b^{-1}c^{-1}b^{-1}cbab^{-1}\\
a^{-1}b_3 = a^{-1}b^{-1}cbaba^{-1}b^{-1}c^{-1}b
\end{array}
$
\vspace{32pt}
}
\bigskip

\noindent Equations on $t_0=\tr(AB)$, $t_1=\tr(AC)$, $t_2=\tr(BC)$:
$$\left\{
\begin{array}{l}
t_0^2t_1^3t_2^2+t_0t_1^2t_2(2-2t_2^2+t_0^2(-2+3t_2^2))+t_0t_2(5-4t_2^2+t_2^4+t_0^4(-1+t_2^2)^2+\\
\ +t_0^2(-4+5t_2^2-2t_2^4))+t_1(1-3t_2^2+t_2^4+t_0^2(-3+7t_2^2-4t_2^4)+t_0^4(1-4t_2^2+3t_2^4)) = 2,\\
t_2+(t_1+t_0t_2)(t_0-t_1t_2-t_0t_2^2) = -2,\\
t_0+(t_1+t_2t_0)(t_2-t_1t_0-t_2t_0^2) = -2.\\
\end{array}
\right.$$
The complex solution:
$$\rho_1 = 2,\ \rho_0 = \rho_2 = \frac{-1\pm\sqrt{-3}}{2}.$$

This example previously appeared as an orbifold which has the sister of the figure-eight knot as a
two-fold cover. As we see, the orbifold is arithmetic with the field of definition 
$\Q[\sqrt{-3}]$ (again by the arithmeticity test from \cite{Bel1}).

\bigskip

\psfig{file=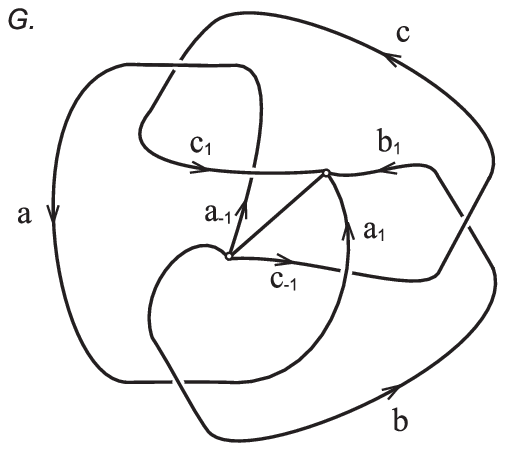, scale=.75}
\hfill
\parbox[b]{0.5\textwidth}{
$\begin{array}{l}
a_1 = bab^{-1}\\
b_1 = cbc^{-1}\\
c_1 = aca^{-1}\\
a_{-1} = c_1^{-1}ac_1\\
c_{-1} = a_1ca_1^{-1}\\
\\
a_1b_1c_1 = bab^{-1}cbc^{-1}aca^{-1}\\
b^{-1}c_{-1}^{-1} = ab^{-1}c^{-1}ba^{-1}b^{-1}\\
a_1b_1 = bab^{-1}cbc^{-1}
\end{array}
$
\vspace{32pt}
}
\bigskip

In our usual notation $t_0=\tr(AB)$, $t_1=\tr(AC)$, $t_2=\tr(BC)$
we obtain 
$$\tr(A_1B_1C_1)=\tr(ABC)(t_0t_1t_2+1) = 0.$$
Since $\tr(ABC)\neq 0$, it gives $t_0t_1t_2+1 = 0,$ and so we have the
following system of equations on the parameters:
$$\left\{
\begin{array}{l}
t_0t_1t_2 = -1,\\
t_0t_2^2+t_1t_2-t_0 = -2,\\
t_0^2t_2+t_1t_0-t_2 = -2.
\end{array}
\right.$$
The complex solution:
$$\rho_0 = \rho_1 = \rho_2 = \frac{1\pm\sqrt{-3}}{2}.$$
Since $\rho_0 = \rho_1 = \rho_2$, the orbifold has an order $3$ symmetry
which cyclically permutes the generators of the group. This symmetry can
be easily seen on the spatial representation of the singular graph but
it is lost on the plane projection. According to \cite{Bel1}, this
orbifold is one of the three arithmetic non-compact generalized triangle
orbifolds with equal parameters (the remaining two orbifolds correspond
to Type~C).

\bigskip

\textbf{\textit{H.}} We have tried more than 20 different structures for this case but
did not succeed in finding any arithmetic examples. The central part structure
implies rather strong restrictions on the orbifold: it should be non-compact
and it is essentially asymmetric. Let us leave as an open problem {\it to find an
arithmetic example for Type~$H$ or to prove that there do not exist such examples.}

\section{ Examples of generalized triangle orbifolds without Property~RE.}

Here we shall consider examples of generalized triangle orbifolds whose
singular sets do not correspond to any of the eight types of singularities
from Theorem~\ref{CT}. 

Our basic example is well known, it is the Picard orbifold. It was first shown in
\cite{Brun} that the Picard group can be generated by only two elements or by
three involutions. Let us recall this presentation.

\bigskip

The singular set of the Picard orbifold is:
\bigskip

\psfig{file=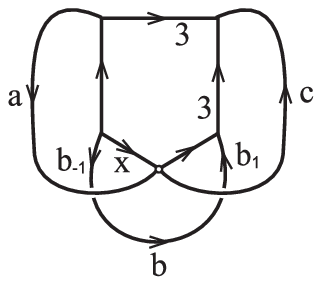, scale=.8}
\hfill
\parbox[b]{0.65\textwidth}{
$\begin{array}{l}
b_1 = cbc^{-1}\\
b_{-1} = aba^{-1}\\
\end{array}
$
\vspace{72pt} }
\bigskip

We shall not write down a matrix representation of the Picard group
because it is well known and can be easily obtained similarly to the examples considered in
the previous section. Our goal is only to show that this group can be generated by 
three involutions, so we shall make all the calculations in the group.

We have:
$$\Gamma =\ <a,b,c,x \mid a^2, b^2, c^2, x^2, (xac)^2, (xabc)^3,(bax)^3, (abax)^2>.$$ 

$(xac)^2 = 1$ implies $acx = xca$;

$(abax)^2 = 1$ implies $xaba = abax$.
\medskip

Now consider the word $(xabc)^3 = 1$: 
$$(xabc)^3 = xabcxabcxabc = xabcxab\ aa\ cxabc = xabc\ abax\ xca\ abc = xabcabacbc;$$ 

$x = abcabacbc.$
\medskip

So $\Gamma =\ <a,b,c \mid a^2, b^2, c^2, (abcabacbc)^2, (abcabacbcac)^2,
(abacb)^3, (acabacbc)^2>.$

\bigskip

Let $\tilde X$ be a $3$-manifold obtained as a complement to the singular set 
of the Picard orbifold in $\S^3$. We have $H_1(\tilde{X})\cong
\Z^4$, so the rank of $\pi_1(\tilde{X})$ is greater than 3, and hence the  
Picard group does not have Property~RE.

The Picard orbifold $\mathcal{O}$ is a non-compact orbifold with a non-rigid cusp. It follows that we can obtain an infinite series 
of compact generalized triangle orbifolds without Property~RE by $(n,0)$-Dehn fillings
on the cusp of $\mathcal{O}$.

\medskip

In \cite{Bel2} the generalized triangle orbifolds without Property~RE were called 
``exceptional''. In fact, the examples considered here have a very special structure
and we expect that there are either no any or very few other examples of 
generalized triangle groups without Property~RE.

\end{document}